\documentclass[12pt,oneside]{amsproc}
\usepackage[T2A]{fontenc}
\usepackage[english]{babel}
\usepackage{amssymb}
\usepackage{amsthm}

\title{Some remarks on spectra of nuclear operators
}
\author{Oleg I. Reinov}
\address{Saint Petersburg State University
}
\email{orein51@mail.ru}

 \thanks{%${ }^\maltese$
}

\begin{document}

  \maketitle

%%%%%%%%%%%%%%%%%%%%%%%%%%%%%%%%%%%

\begin{abstract}
It was shown by M. I. Zelikin (2007) that the spectrum of a nuclear operator in
a separable Hilbert space is central-symmetric iff the spectral traces of all
odd powers of the operator equal zero. The criterium can not be extended to the case of general Banach spaces:
It follows from Grothendieck-Enflo results that
there exists a nuclear operator $U$ in the space $l_1$ with the property that
 $\operatorname{trace}\, U=1$ and $U^2=0.$
B. Mityagin (2016) has generalized Zelikin's criterium to the case of compact operators
(in Banach spaces) some of which powers are nuclear.
We give sharp generalizations of Zelikin's theorem (to the cases of subspaces of quotients of $L_p$-spaces)
and of Mityagin's result (for the case where the operators
are not necessarily compact). 
\end{abstract}
 \medskip

%%%%%%%%%%%%%%%%%%%%%

{\bf 1.\, Introduction}.\
It was shown by M. I. Zelikin in [16] that the spectrum of a nuclear operator in
a separable Hilbert space is central-symmetric if and only if the spectral traces of all
odd powers of the operator equal zero. Recall that the spectrum of every nuclear operator
in a Hilbert space consists of non-zero
eigenvalues of finite algebraic
 multiplicity, which have no limit point except possibly zero, and maybe zero.
This system of all eigenvalues (written according to their multiplicities) is
absolutely summable, and the spectral trace of any nuclear operator is, by definition,
the sum of all its eigenvalues (taken according to their multiplicities).

The space of nuclear operators in a Hilbert space may be defined as the space
of all trace-class operators (see [14, с. 77]);  % [3, Theorem 8.1];   %!!! check eng pages
in this case we speak about the "nuclear trace" of an operator).
Trace-class operators in a Hilbert space can be considered also as the elements of the completion
of the tensor product of the Hilbert space and its Banach dual with respect to the greatest
crossnorm on this tensor product [14, с. 119].
The well-known Lidski\v{\i} theorem [5] %(see also [3, Theorem 8.4])   %!!! eng pages and V.\,B.~Lidski\v{\i}
says that the nuclear trace of any nuclear operator in a Hilbert space
(or, what is the same, of the corresponding tensor element) coincides with
its spectral trace. Thus, Zelikin's theorem [16] can be reformulated in the following way:
the spectrum of a nuclear operator in a separable Hilbert space is central-symmetric if and only if
the nuclear traces of all odd powers of the corresponding tensor element are zero.

One of the aim of our notes is to give an exact generalization of this result to the case
of tensor elements of so-called $s$-projective tensor products of subspaces of
quotients of  $L_p(\mu)$-spaces. In particular, we get as a consequence Zelikin's theorem
(taking $p=2).$

Another problem which is under consideration in our notes is concentrated around
the so-called $\mathbb Z_d$-symmetry of the spectra of the linear operators.
The notion of $\mathbb Z_d$-symmetry of the spectra was introduced by B. S. Mityagin
in a preprint [6] %!!!REF
and in his paper [7].  %!!!REF
He is interested there in a generalization of the result from [16] %!!!REF
in two directions: to extend Zelikin's theorem to the case of general Banach
spaces and to change the property of a compact operator to have
central-symmetric spectrum %to the property of the one
to have $\mathbb Z_d$-symmetric spectrum. Roughly speaking,
$\mathbb Z_d$-symmetry of a spectrum of a compact operator $T$ means
that for any non-zero eigenvalue $\lambda$ of $T$ the spectrum contains also
as eigenvalues of the same algebraic multiplicities all "$d$-shifted"\, numbers
$t\lambda$ for $t\in \sqrt[d]{1}.$
B. S. Mityagin has obtained a very nice result, showing that the spectrum of a compact operator $T$ in an
arbitrary Banach space, some power $T^m$ of which is nuclear,
 is $\mathbb Z_d$-symmetric if and only if for all large enough integers of type $kd+r$
$(0< r < d)$  the nuclear traces of $T^{kd+r}$ are zero.
We present some thoughts around this theorem, giving, in particular, a short
(but not so elementary as in [6, 7]) proof for the case where the operator
is not necessarily compact. Let us mention, however, that the proof from [6, 7] %!!!REF
can be adapted for this situation too.

Some words about the content of the paper.

In Section 2, we introduce some notation, definitions and terminology in connection
with so-called $s$-projective tensor products, $s$-nuclear operators and the
approximation properties of order $s,$\, $s\in (0, 1].$ We formulate here two
auxiliary assertions from the author's paper [13]; they give us generalized
Grothendieck-Lidski\v{\i} trace formulas which will be useful
in the next section.

Section 3 contains an exact generalization of Zelikin's theorem. In this
section we present a criterion for the spectra of $s$-nuclear operators
in subspaces of quotients of $L_p$-spaces to be central-symmetric.

Results of Section 4 show that the criterion of the central symmetry, obtained
in the previous section, is optimal. In particular, we present here (Theorem 2)
sharp examples of $s$-nuclear operators $T$ in the spaces $l_p,$\, $1\le p\le+\infty,$\,
$p\neq2,$ for which $\operatorname{trace}\, T=1$ and $T^2=0.$

Finally, Section 5 is devoted to the study of Mityagin's $\mathbb Z_d$-symmetry
of the spectra of linear operators.
Our aim here is to give a short (but using the Fredholm Theory) proof
of Mityagin's theorem [6, 7] for arbitrary linear continuous (Riesz) operators.
Firstly, we consider a $\mathbb Z_2$ situation
(central symmetry) to clarify an idea which is to be used then in the general case.
 We finish the paper with a short proof of the theorem from [6-7] for
 continuous (not necessarily compact) operators and with some simple examples
 of applications.

   \bigskip

      %%%%%
      {\bf 2.\, Preliminaries}.\
We denote by $X, Y, \dots$ Banach spaces,
$L(X, Y)$ is a Banach space of all linear continuous operators
from $X$  to $Y;$ $L(X):= L(X, X).$
For a Banach dual to a space $X$ we use the notation $X^*.$
  If $x\in X$ and $x'\in X^*,$ then $\langle x',x\rangle$ denotes the value $x'(x).$

      By   $X^*\widehat\otimes X$
   we denote the projective tensor product of the spaces $X^*$ and  $X$ [2] (see also [11, 12]).    %% !!!
  It is a completion of the algebraic tensor product
   $X^*\otimes X$  (considered as a linear space of all finite rank continuous
   operators $w$ in $X)$  with respect to the norm
  $$
  ||w||_{\land}:= \inf \{\left(\sum_{k=1}^N ||x'_k||\, ||x_k||\right):\ w=\sum_{k=1}^N x'_k\otimes x_k\}.
  $$
   Every element
   $u$ of the projective tensor product $X^*\widehat\otimes X$
     can be represented in the form
    $$
    u=\sum_i \lambda_i x'_i\otimes x_i,
    $$
where $(\lambda_i)\in l_1$ и $||x'_i||\le1,$ $||x_i||\le1$ [2].  %!!!

   More generally, if $0<s\le 1,$ then
   $X^*\widehat\otimes_s X$ is a subspace of the projective tensor product,
   consisting of the tensor elements $u, u\in X^*\widehat\otimes X,$
   which admit representations of the form
   $u=\sum_{k=1}^\infty x'_k\otimes x_k,$
   where $(x'_k)\subset X^*, (x_k)\subset X$ and
   $\sum_{k=1}^\infty ||x'_k||^s\, ||x_k||^s <\infty$ [2, 11, 12].
   Thus,
   $X^*\widehat\otimes X = X^*\widehat\otimes_1 X.$

On the linear space $X^*\otimes X,$ a linear functional "trace" is defined by a natural way.
It is continuous on the normed space %!!! " " OK? in PDF - ok.
$(X^*\otimes X, ||\cdot||_{\land})$ and has the unique continuous extension
    to the space $X^*\widehat\otimes X,$ which we denote by $\operatorname{trace}.$

Every tensor element $u, u\in X^*\widehat\otimes X,$ of the form
    $u=\sum_{k=1}^\infty x'_k\otimes x_k$
 generates naturally an operator
    $\widetilde u: X\to X,$ $\widetilde u(x):= \sum_{k=1}^\infty \langle x'_k, x\rangle\,  x_k$ for
    $x\in X.$
 This defines a natural mapping $j_1: X^*\widehat\otimes X \to L(X).$ The operators,
    lying in the image of this map, are called nuclear [2],  [8].  %!!!
   More generally, if $0<s\le1,$  $u=\sum_{k=1}^\infty x'_k\otimes x_k$ and
    $\sum_{k=1}^\infty ||x'_k||^s\, ||x_k||^s <\infty,$ then the corresponding operator $\widetilde u$
    is called  $s$-nuclear [11, 12]. By $j_s$ we denote a natural map from
    $X^*\widehat\otimes_s X$ to $L(X).$
   We say that a space $X$ has the approximation property of order
    $s,$\, $0<s\le1$ (the $AP_s),$ if the canonical mapping
    $j_s$ is one-to-one [11, 12]. Note that the $AP_1$ is exactly the approximation
    property $AP$ of A. Grothendieck [2], [8].   %!!!
   Classical spaces, such as $L_p(\mu)$  and $C(K),$ have the approximation property.
  If a space $X$ has the $AP_s,$ then we can identify the tensor product
    $X^*\widehat\otimes_s X$  with the space $N_s(X)$  of all
    $s$-nuclear operators in $X$
   (i.e. with the image of this tensor product under the map $j_s).$ In this case
   for every operator $T\in N_s(X)= X^*\widehat\otimes_s X$  the functional $\operatorname{trace}\, T$
   is well defined and called
   the  {\it nuclear trace}\, of the operator $T.$

  It is clear that if a Banach space has the approximation property, then
  it has all the properties $AP_s,$ $s\in (0,1].$
    Every Banach space has the property $AP_{2/3}$ (A. Grothendiek [2], see also [10]).  %!!!
  Since each Banach space is a subspace of an $L_\infty(\mu)$-space,
    the following fact (to be used below) is a generalization of the mentioned
    result of A. Grothendieck:

    \smallskip
 %%%%%%%%%%%%%%%%%%%%%%%
{\bf Lemma 1} [13, Corollary 10].\,    %!!!
  Let $s\in (0,1],$ $p\in [1,\infty]$ and $1/s=1+|1/p-1/2|.$
    If a Banach space $Y$ is isomorphic to a subspace of a quotient
(or to a quotient of a subspace) of some
        $L_p(\mu)$-space, then it has the $AP_s.$

       \vskip 0.3cm

 Thus, for such spaces we have an equality $Y^*\widehat\otimes_s Y= N_s(Y)$ and
  the nuclear trace of any operator $T\in N_s(Y)$
  is well defined.

  We will need also the following auxiliary assertion
  (the first part of which is a consequence of the previous lemma).

  {\bf Lemma 2} [13, Theorem 1].\,     %!!!
Let $Y$ be a subspace of a quotient
(or a quotient of a subspace) of some
        $L_p(\mu)$-space, $1\le p\le\infty.$
If  $T\in N_s(Y),$  where  $1/s=1+|1/2-1/p|,$
then

1.\, the nuclear trace of the operator $T$ is well defined,

2.\, $\sum_{n=1}^\infty |\lambda_n(T)|<\infty,$ where                     %!!! notation before?
$\{\lambda_n(T)\}$ is the system of all eigenvalues of the operator $T$
(written according to their algebraic multiplicities), and
$$
 \operatorname{trace}\, T= \sum_{n=1}^\infty \lambda_n(T).
$$

\vskip 0.3cm

      %%%

Following  [16], we say that a spectrum of a compact operator in a Banach space
is
{\it central-symmetric}, if
for every its eigenvalue
$\lambda$ the number $-\lambda$ is also its eigenvalue and of the same algebraic multiplicity.
 We shall use the same terminology in the case of operators
all non-zero spectral values of which %$\lambda(\widetilde{v})\in \operatorname{sp}\, (\widetilde{v})$
are eigenvalues of finite
 multiplicity and have no limit point except possibly zero; the corresponding
 eigenvalue sequence for such an operator $T$ will be denoted by $\operatorname{sp}\, (T);$
 thus it is an unordered sequence of all eigenvalues of $T$ taken according to their multiplicities.
  \bigskip

                %%%%%%%%%%%%%%
                {\bf 3.\, On central symmetry}.\
Let us note firstly that the theorem of Zelikin (in the form as it was formulated in [16])
can not be extended to the case of general Banach spaces, %!!!
even if the spaces have the Grothendieck approximation property.

{\bf Example 1}.\
Let  $U$ be a nuclear operator in the space $l_1$, constructed in [8, Proposition 10.4.8].
This operator has the property that $\operatorname{trace}\, U=1$ and
$U^2=0.$ Evidently, the spectrum of this operator is $\{0\}.$ Let us note that the operator
is not only nuclear, but also belongs to the space $N_s(l_1)$ for all $s\in (2/3, 1].$
It is not possible to present such an example in the case of $2/3$-nuclear operators
(see Corollary 3 below).
Note also that, however, the traces of all operators
$U^m, m=2, 3, \dots,$ (in particular, $U^{2n-1})$ are equal to zero.
\smallskip

{\bf Remark 1}.\
For every nuclear operator $T: X\to X$  and for any natural number $n>1,$ the nuclear trace
$\operatorname{trace}\, T^n$ is well defined (see [2, Chap. II, Cor. 2, p. 16]) and
equal  the sum of all its eigenvalues (according to their multiplicity)
[2, Chap. II,  Cor. 1, p. 15]. Therefore, if the spectrum
of a nuclear operator $T: X\to X$ is central-symmetric, then for each odd
$m=3, 5, 7, \dots$ the nuclear trace of the operator
$T^m$ is equal to zero. This follows from the fact that
the eigenvalue sequences of $T$ and $T^m$ can be arranged in such a way that
$\{\lambda_n(T)^m\}= \{\lambda_n(T^m)\}$ (see, e.g., [9, 3.2.24, p. 147]).
\smallskip

        Let us formulate and prove the central result of this section.
                \smallskip

{\bf Theorem 1}.\ {\it
 Let  $Y$ be a subspace of a quotient (or a quotient of a subspace)
of an $L_p$-space, $1\le p\le\infty,$
and $u\in Y^*\widehat\otimes_s Y,$  where $1/s=1+|1/2-1/p|,$
The spectrum of the operator $\widetilde u$
is central-symmetric if and only if
$\operatorname{trace}\, u^{2n - 1} = 0, n =1,2,\dots.$
}
\smallskip

{\it Proof}.\
If the spectrum of $\widetilde u$ is central-symmetric, then, by Lemma 2,
$\operatorname{trace}\, u= \operatorname{trace}\, \widetilde u= \sum_{n=1}^\infty \lambda_n(T)=0;$ also, by Remark 1,
  $\operatorname{trace}\, u^m= \sum_{n=1}^\infty \lambda_n(T^m)=0$ for $m=3, 5, \dots.$

 To prove the inverse, we need some information from the Fredholm Theory.
Let $u$ be an element of the projective tensor product $X^*\widehat\otimes X,$
where $X$ is an arbitrary Banach space.
 Recall that the Fredholm determinant $\operatorname{det}\, (1-zu)$ of $u$
 (see [2, Chap. II,  p. 13] or [3, 8, 9])
  is an entire function
$$
 \operatorname{det}\, (1-zu)= 1-z\, \operatorname{trace}\, u+ \dots + (-1)^n z^n\alpha_n(u) +  \dots,
$$
all zeros of which are exactly (according to their multiplicities) the inverses of nonzero eigenvalues %$(\mu_k)$
of  the operator $\widetilde{u},$ associated with the tensor element $u.$
By [2, Chap. II, Cor. 2, pp. 17-18],
this entire function is of the form
$$
  \operatorname{det}\, (1-zu) =
  e^{-z\operatorname{trace}\, u} \prod_{i=1}^\infty (1- z z_i)\, e^{zz_i},
$$
where
$\{z_i=\lambda_i(\widetilde u)\}$ is a system of all eigenvalues of the operator $\widetilde{u}$
(written according to  their algebraic multiplicities).
Also, there exists a $\delta>0$ such that for all $z, |z|\le\delta,$ we have
%$$
%    \operatorname{det}\, (1-zu) =
%  \exp {(-\sum_{n=1}^\infty \frac1n\, z^n \operatorname{trace}\, u^n)}
%  $$

\begin{equation}
     \operatorname{det}\, (1-zu) =
  \exp {(-\sum_{n=1}^\infty \frac1n\, z^n \operatorname{trace}\, u^n)}         \label{(1)}
\end{equation}
(see [3, p. 350]; cf. [1, Theor. I.3.3, p. 10]).      %!!!Ref

Now,
let $u\in Y^*\widehat\otimes_s Y$ be as in the formulation of our theorem and
suppose that $\operatorname{trace}\, u^{2n - 1} = 0, n =1,2,\dots.$
 By (\ref{(1)}), we get: for a neighborhood $U= U(0)$ of zero in $\mathbb C,$
 $\operatorname{det}\, (1-zu)= \operatorname{det}\, (1+zu)$ for $z\in U.$ 
 Therefore, the entire function $\operatorname{det}\, (1-zu)$
 is even.
  By definition of $\operatorname{det}\, (1-zu),$ the sequence
of zeros of this function is exactly the sequence of inverses of nonzero
eigenvalues of $\widetilde u.$
Hence, the spectrum of $\widetilde u$ is central-symmetric.
  \smallskip

    Since under the conditions of Theorem 1 the space $Y$ has the $AP_s,$
   the tensor product $Y^*\widehat\otimes_s Y$ can be identified naturally
   with the space of all $s$-nuclear operators in $Y.$ Hence, the statement of Theorem 1
   may be reformulated in the following way:
\smallskip

   {\bf Corollary 1}. {\it
      Let $s\in [2/3,1], p\in[2,\infty], 1/s=1+|1/2-1/p|,$
   $Y$ be a subspace of a quotient (or a quotient of a subspace)
of an $L_p$-space,
$T$ be an $s$-nuclear operator in $Y.$
The spectrum of  $T$ is central-symmetric if and only if
$\operatorname{trace}\,  T^{2n - 1} = 0, \, n \in \mathbf N.$
   }
\smallskip

      {\bf Corollary 2}\, [16]. {\it   %!!!
The spectrum of a nuclear operator $T,$ acting on a Hilbert space,
is central-symmetric if and only if
$\operatorname{trace}\,  T^{2n - 1} = 0, \, n \in \mathbf N.$
   }

   For a proof, it is enough to apply Theorem 1 for the case $p=2.$
\smallskip

{\bf Corollary 3}. {\it
The spectrum of a 2/3-nuclear operator $T,$ acting on an arbitrary Banach space,
is central-symmetric if and only if
$\operatorname{trace}\,  T^{2n - 1} = 0, \, n \in \mathbf N.$
   }

   For a proof, it is enough to apply Theorem 1 for the case $p=\infty,$
    taking into account the fact that every Banach space is isometric to
   a subspace of an $L_\infty(\mu)$-space.
 \smallskip

   In connection with Corollary 3, let us pay attention again to the nuclear operator
   from Example 1.
    \medskip

    %%%%%%%%%%%%%%
                {\bf 4.\,  Sharpness of results of Section 3}.\
            %%%%%%%%%%  $1/s=1+|1/2-1/p|,$
   Now we will show that the statement of Theorem 1 is sharp and that
   the exponent $s$ is optimal if $p$ is fixed   (if of course $p\neq2,$ i.e. $s\neq1).$

Consider the case $2< p\le\infty.$ In this case $1/s=1+|1/2-1/p|= 3/2-1/p.$
   In a paper of the author [12, Example 2] %!!!
 the following result was obtained (see a proof in [12]):  %!!!

 ${(\star)}$\,
    {\it
 Let $r\in [2/3,1), p\in(2,\infty], 1/r=3/2-1/p.$
There exist a subspace
 $Y_p$ of the space $l_p$
$(c_0$ if $p=\infty)$
and a tensor element
$w_p\in Y_p^*\widehat\otimes_1 Y_p$ such that
$w_p\in Y_p^*\widehat\otimes_s Y_p$ for every $s>r,$
$\operatorname{trace}\, w_p=1, \tilde w_p=0$  and the space $Y_p$ (as well as $Y_p^*)$
has the $AP_r$
$($but evidently does not have the $AP_s$ if $1\ge s>r).$
Moreover, this element admits a nuclear representation of the form
$$
  w_p=\sum_{k=1}^\infty \mu_k\, x'_k \otimes x_k,\  \text{ where }\ ||x'_k||=||x_k||=1,\,
  \sum_{k=1}^\infty |\mu_k|^s<\infty\ \forall \, s>r.
$$
}
%\smallskip

Evidently, we have for a tensor element $u:=w_p$ from the assertion ${(\star)}:$
$\operatorname{trace}\, u=1$ and the spectrum of the operator $\widetilde u$ equals $\{0\}.$

The case where  $2< p\le\infty$ can be considered analogously (with an application of the assertion ${(\star)}$
to a  "transposed" tensor element
$w_p^t\in Y_p \widehat\otimes_1 Y_p^*.)$

   As was noted above (Example 1), there exists a nuclear operator $U$ in $l_1$ such that
$U^2=0$ and  $\operatorname{trace}\, U=1.$ The following theorem is an essential
generalization of this result and gives us the sharpness of the statement of Corollary 1
(even in the case where $Y=l_p).$
\smallskip

       {\bf Theorem 2}.\,  {\it
       Let  $p\in [1,\infty],  p \neq2,$ $1/r= 1+|1/2-1/p|.$
There exists a nuclear operator $V$ in $l_p$ (in $c_0$ if $p=\infty)$
such that

1)\,
$V \in N_s(l_p)$ for each $s\in (r, 1];$

2)\,
$V\notin N_r(l_p);$

3)\,
$\operatorname{trace}\, V=1$ and $V^2=0.$
}
\smallskip

{\it Proof}.\,
 Suppose that $p>2.$
Consider the tensor element $w:= w_p$ from the assertion $(\star)$ and its representation
$  w=\sum_{k=1}^\infty \mu_k\, x'_k \otimes x_k,$ where
$||x'_k||=||x_k||=1$ и $\sum_{k=1}^\infty |\mu_k|^s<\infty$ for each $s, s>r.$
 Let $l: Y:= Y_p\to l_p$ be the identity inclusion. Let $y'_k$ be
 an extension of the functional $x'_k$\, $(k=1, 2, \dots)$ from the subspace
 $Y$ to the whole space $l_p$ with the same norm and set
 $v:= \sum_{k=1}^\infty \mu_k\, y'_k \otimes l(x_k).$
 Then $v\in l_{p'} \widehat\otimes_s l_p$ \, $(1/p+1/p'=1)$\,  for each $s\in (r, 1],$
 $\operatorname{trace}\, v= \sum \mu_k\, \langle y'_k, l(x_k)\rangle=1$ and
 $\widetilde v(l_p)\subset l(Y)\subset l_p.$
 On the other hand, we have a diagram:
 $$
  Y \overset{l}\to l_p \overset{\widetilde v_0}\to Y \overset{l}\to l_p \overset{\widetilde v_0}\to Y \overset{\widetilde v_0}\to l_p,
 $$
where $\widetilde v_0$ is an operator generated by $\widetilde v,$  $\widetilde v = l\widetilde v_0$  and
$\widetilde v_0 l= \widetilde w=0.$
Put $V:= \widetilde v.$ Clearly,  $\operatorname{trace}\, V=1$ and the spectrum
$\operatorname{sp}\, \,V^2=\{0\}.$ Let us note that $V\notin N_r(l_p)$ (by Lemma 2).
  If $p\in [1,2),$ then it is enough to consider the adjoint operator.
 \smallskip

 It follows from Theorem 2 that the assertion of Corollary 1 is optimal already in the case
of the space $Y=l_p$ (which, by the way, has the Grothendieck approximation property).
    \bigskip

    %%%%%%%%%%%%%%
                {\bf 5.\,  Generalizations: around Mityagin's theorem}.\
Recall that if $T\in L(X)$ and, for some $m\in \mathbb N,$ $T^m$ is a Riesz operator
(see, e.g., [4, p. 943] for a definition), then $T$ is a Riesz operator too
(see, e.g., [9, 3.2.24, p. 147]). In particular, if $T^m$ is compact, then
all non-zero spectral values $\lambda(T)\in \operatorname{sp}\, (T)$ are eigenvalues of finite
(algebraic) multiplicity and have no limit point except possibly zero.
Also, in this case
the eigenvalue sequences of $T$ and $T^m$ can be arranged in such a way that
$\{\lambda_n(T)^m\}= \{\lambda_n(T^m)\}$ (see [9, 3.2.24, p. 147]).
   Recall that in this case  we denote by $\operatorname{sp}\, (T)$ 
    (resp., by $\operatorname{sp}\, (T^m)$ the sequence
$\{\lambda_n(T)^m\}$  (resp., $\{\lambda_n(T^m)\}).$
        %  \smallskip

 We are going to present a short proof of the theorem
of B. Mityagin from [6, 7].   %!!!REF  and below
To clarify our idea of the proof, let us consider firstly the simplest case
where $d=2.$
           \smallskip

{\bf Theorem 3.}\ {\it
Let $X$ be a Banach space and $T\in L(X).$ Suppose that some power of $T$
is nuclear. The spectrum of $T$ is central-symmetric if and only if there is an integer $K\ge0$ such that
 for every $l> K$ the value $\operatorname{trace}\, T^{l}$ is well defined and
 $\operatorname{trace}\, T^{2l+1}=0$ for all $l> K.$
 }
     \smallskip

     {\it Proof}.\
    Suppose that $T\in L(X)$ and there is an $M\in \mathbb N$ so that $T^M\in N(X).$
Fix an odd $N_0, N_0>M,$ with the property that $T^{N_0}\in N_{2/3}(X)$
(it is possible since a product of three nuclear operators is 2/3-nuclear)
and $\operatorname{trace}\, T^{N_0+2k}=0$ for all $k= 0, 1, 2, \dots.$
By Corollary 3, the spectra of all $T^{N_0+2k}$ are central-symmetric
(since, e.g., $\operatorname{trace}\, T^{N_0}= \operatorname{trace}\, (T^{N_0})^3 = \operatorname{trace}\, (T^{N_0})^5= \dots = 0$
by assumption). Assume that the spectrum of $T$ is not central-symmetric.
Then there exists an eigenvalue $\lambda_0\in \operatorname{sp}\, (T)$ such that $-\lambda_0\notin \operatorname{sp}\, (T).$

Now, $\lambda_0^{N_0}\in \operatorname{sp}\, (T^{N_0}),$ so $-\lambda_0^{N_0}\in \operatorname{sp}\, (T^{N_0}).$
Hence,
there exist $\mu_{N_0}\in \operatorname{sp}\, (T)$ and $\theta_{N_0}$ so that $|\theta_{N_0}|=1,$
$\mu_{N_0}^{N_0}=-\lambda_0^{N_0}$ and $\mu_{N_0}=\theta_{N_0} \lambda_0,$\, $\theta_{N_0}\neq -1.$
Analogously,
     $\lambda_0^{N_0+2}\in \operatorname{sp}\, (T^{N_0+2}),$ so $-\lambda_0^{N_0+2}\in \operatorname{sp}\, (T^{N_0+2}).$
Hence,
there exist $\mu_{N_0+2}\in \operatorname{sp}\, (T)$ and $\theta_{N_0+2}$ so that $|\theta_{N_0+2}|=1,$
$\mu_{N_0+2}^{N_0+2}=-\lambda_0^{N_0+2}$ and $\mu_{N_0+2}=\theta_{N_0+2} \lambda_0,$
\, $\theta_{N_0+2}\neq -1$ etc.
By induction we get the sequences   $\{\mu_{N_0+2k}\}_{k=0}^\infty$ and
$\{\theta_{N_0+2k}\}_{k=0}^\infty$ with the properties that
  $\mu_{N_0+2k}\in \operatorname{sp}\, (T),$\, $|\theta_{N_0+2k}|=1,$
$\mu_{N_0+2k}^{N_0+2k}=-\lambda_0^{N_0+2k}$ and $\mu_{N_0+2k}=\theta_{N_0+2k} \lambda_0,$
\, $\theta_{N_0+2k}\neq -1.$
Since $\mu_{N_0+2k}\in \operatorname{sp}\, (T)$\, and $|\mu_{N_0+2k}|=|\lambda_0|>0,$
the sequence $\{\mu_{N_0+2k}\}$ is finite as a set, i.e.,
we have that $\mu_{N_0+2K}=\mu_{N_0+2K+2}=\dots$ for some $K>1.$
It follows that $\theta_{N_0+2K}=\theta_{N_0+2K+2}=\dots.$
But $\theta_{N_0+2k}^{N_0+2k}=-1$ for all $k.$ Thus
$\theta_{N_0+2K}^l=-1$ for every odd $l\ge N_0+2K.$ Therefore,
  $\theta_{N_0+2K}=-1.$ Contradiction.
        \medskip

Now we are going to consider a general case of a notion of $\mathbb Z_d$-symmetry
of a spectra,
introduced and investigated by B. Mityagin in [6, 7].  %!!!REF
Let $T$ be an operator in $X,$
all non-zero spectral values of which %$\lambda(\widetilde{v})\in \operatorname{sp}\, (\widetilde{v})$
are eigenvalues of finite
 multiplicity and have no limit point except possibly zero. Recall that we denote by
 $\operatorname{sp}\, (T)$ the corresponding unordered
 eigenvalue sequence for  $T$ (possibly, including zero).
For a fixed $d=2, 3, \dots$ and for the operator $T,$ the spectrum of $T$
 is called $\mathbb Z_d$-symmetric,
if $\lambda\in \operatorname{sp}\, (T)$ implies $t\lambda\in \operatorname{sp}\, (T)$ for every $t\in\sqrt[d]{1}.$

   Let $r\in (0,\infty],$ \, $\mathbb D:= \{z\in \mathbb C:\  |z|<r\},$
   $f: \, \mathbb D\to \mathbb C.$  and $d\in \mathbb N\setminus\{1\}.$
   We say that $f$ is $d$-even if $f(tz)=f(z)$ for every $t\in \sqrt[d]{1}.$
   \smallskip

   {\bf Lemma 3.}\                                     %!!!
 Let $\Phi(X)$ be a linear subspace of $X^*\widehat\otimes X$ of spectral type $l_1,$
 i.e., for every $v\in \Phi(X)$ the series
% all non-zero spectral values $\lambda(\widetilde{v})\in \operatorname{sp}\, (\widetilde{v})$ are eigenvalues of finite
 % multiplicity, have no limit point except possibly zero and
  $\sum_{\lambda\in\operatorname{sp}\, (\widetilde{v})} |\lambda|$ is convergent. Let $d\in\mathbb N, d>1.$
 If $u\in\Phi(X),$ then the Fredholm determinant $\operatorname{det}\, (1-zu)$ is $d$-even if and only if
 the eigenvalue sequence of $\widetilde u$ is $\mathbb Z_d$-symmetric if and only if
 $\operatorname{trace}\, u^{kd+r}=0$ for all $k=0, 1, 2, \dots$ and $r= 1, 2, \dots, d-1.$
    \smallskip

    {\it Proof}.\,
  If the function $\operatorname{det}\, (1-zu)$ is $d$-even, then
 the eigenvalue sequence of $\widetilde u$ is $\mathbb Z_d$-symmetric,
 since this sequence coincides with the  sequence of inverses of
  zeros
 of $\operatorname{det}\, (1-zu)$ (according to  their multiplicities).

 If the eigenvalue sequence of $\widetilde u$ is $\mathbb Z_d$-symmetric, then
 $\operatorname{trace}\, u= \sum_{\lambda\in\operatorname{sp}\, (\widetilde{u})} \lambda=0$ (since $\Phi(X)$ is of spectral type $l_1$
 and $\sum_{t\in\sqrt[d]{1}} t=0).$ Also, by the same reason
    $\operatorname{trace}\, u^{kd+r}=0$ for all $k=0, 1, 2, \dots$ and $r= 1, 2, \dots, d-1,$
 since the spectrum of $\widetilde{u}^l$ is absolutely summable for every $l\ge2$
 and we may assume that $\{\lambda_m(\widetilde{u}^l)\}= \{\lambda_m(\widetilde{u})^l\}$ (for every fixed $l)$
 [9, 3.2.24, p. 147].

 Now, let $\operatorname{trace}\, u^{kd+r}=0$ for all $k=0, 1, 2, \dots$ and $r= 1, 2, \dots, d-1.$
 By (\ref{(1)}),
  $\operatorname{det}\, (1-zu) =
  \exp {(-\sum_{m=1}^\infty \frac1{md}\, z^{md} \operatorname{trace}\, u^{md})}$
   in a neighborhood $U$ of zero.
Therefore, this function is $d$-even in the neighborhood $U.$
By the uniqueness theorem $\operatorname{det}\, (1-zu)$ is $d$-even in $\mathbb C.$
     \smallskip

  {\bf Corollary 4.}\   %!!!
For any Banach space $X$ and for every $u\in X^*\widehat\otimes_{2/3} X$ the conclusion of
Lemma 3 holds.                                                                      %!!!
            \smallskip

 {\bf Corollary 5.}\   
 Let $Y$ be a subspace of a quotient of an $L_p$-space, $1\le p\le\infty.$
For any $u\in Y^*\widehat\otimes_s Y,$  where $1/s=1+|1/2-1/p|,$
 the conclusion of Lemma 3 holds.
\medskip

Now we are ready to present a short proof of the theorem
of B. Mityagin [6, 7].   %!!!REF  and below
Note that the theorem in [6, 7] is formulated and proved for compact operators, but
the proof from [6, 7] can be easily adapted for the general case of linear operators.
\smallskip

{\bf Theorem 4.}\
Let $X$ be a Banach space and $T\in L(X).$ Suppose that some power of $T$
is nuclear. The spectrum of $T$ is $\mathbb Z_d$-symmetric if and only if there is an integer $K\ge0$ such that
 for every $l> Kd$ the value $\operatorname{trace}\, T^{l}$ is well defined and
 $\operatorname{trace}\, T^{kd+r}=0$ for all $k=K, K+1, K+2, \dots$ and $r= 1, 2, \dots, d-1.$
     \smallskip

     {\it Proof}.\
 Fix $N_0\in\mathbb N$ such that $T^{N_0}$ is $2/3$-nuclear (it is possible
 by a composition theorem from
[2, Chap, II,  Theor. 3, p. 10]). Note that, by A. Grothendieck,
the trace of $T^l$ is well defined for all $l\ge N_0.$

Suppose that the spectrum of $T$ is $\mathbb Z_d$-symmetric.
Take an integer $l:= kd+r\ge N_0$ with $0<r<d.$
Since the spectrum of ${T}^l$ is absolutely summable,
$\operatorname{trace}\, T^l= \sum_{\lambda\in \operatorname{sp}\, (T^l)} \lambda$
 and we may assume that $\{\lambda_m(T^l)\}= \{\lambda_m(T)^l\},$ we get that
 $\operatorname{trace}\, T^{kd+r}=0.$

  In proving the converse, we may (and do) assume that $Kd>N_0.$
  Consider an infinite increasing sequence $\{p_m\}$ of prime numbers
  with $p_1> (K+1)d.$ Assuming
that
 $\operatorname{trace}\, T^{kd+r}=0$ if $k=K, K+1, K+2, \dots$ and $r= 1, 2, \dots, d-1,$
  for a fixed $p_m$  we get from Lemma 3 (more precisely, from
  Corollary 4) that the function
  $\operatorname{det}\, (1-zT^{p_m})$ is $d$-even.
    %and for all small enough $z$ we have (by (\ref{(1)})
%  $$
 %    \operatorname{det}\, (1-zT^{p_m}) =
  %\exp {(-\sum_{n=1}^\infty \frac1n\, z^n \operatorname{trace}\, (T^{p_m})^n)} =
   %  \exp {(-\sum_{j=1}^\infty \frac1{dj}\, z^{dj} \operatorname{trace}\, T^{p_m dj})}
  %$$
 Suppose that  the spectrum of $T$ is not $\mathbb Z_d$-symmetric.
 Then there exist an eigenvalue $\lambda_0\in \operatorname{sp}\, (T)$ and a root $\theta\in\sqrt[d]{1}$    
 so that $\theta \lambda_0\notin \operatorname{sp}\, (T).$
 On the other hand, again by Lemma 3, the spectrum of $T^{p_m}$ is  $\mathbb Z_d$-symmetric.
 Since $\lambda_0^{p_m}\in \operatorname{sp}\, (T^{p_m}),$ there exists $\mu_m\in \operatorname{sp}\, (T)$ such that
 $\mu_m^{p_m}=\theta^{p_m} \lambda_0^{p_m}\in \operatorname{sp}\, (T^{p_m});$ hence, $\mu_m = \theta_m \lambda_0$
 for some $\theta_m$ with $|\theta_m|=1.$
 But $|\lambda_0|>0.$ Therefore, the set $\{\mu_m\}$ is finite and it follows that there is an integer $M>1$
 such that $\theta_M=\theta_{M+1}=\theta_{M+2}=\dots.$
 Hence, $\theta^{p_m}= \theta_M^{p_m}$
 for all $m\ge M.$ Thus, $\theta_M=\theta.$ Contradiction.
      \smallskip

 Let us give some examples in which we can apply Theorem 4, but the main result of [6, 7]
 does not work.
 \smallskip

 {\bf Example 2.}\
 Let $\Pi_p$ be the ideal of absolutely $p$-summing operators $(p\in [1,\infty);$
see [8] for a definition and related facts). Then for some $n$ one has
 $\Pi_p^n\subset N.$ In particular, $\Pi_2^2(C[0,1])\subset N(C[0,1]),$ but
 not every absolutely 2-summing operator in $C[0,1]$ is compact. Another interesting example:
 $\Pi_3^3$ is of spectral type $l_1$ [15]. %!!!REF
 We do not know (maybe, it is not known to everybody), whether the finite rank operators
 are dense in this ideal. However, Theorem 4 may be applied. Moreover, it can be seen that,
for example, the spectrum of an operator $T$ from $\Pi_3^3$ is central-symmetric if and only if
{\it the spectral traces}\, of the operators $T^{2k-1}$ are zero for all $k>0.$
    \bigskip

{\bf Acknowledgement.}\,
     I would like to express my gratitude to B. S. Mityagin
for a two-days helpful discussion on the topic of this paper and on the
problems related to the notion of nuclear operators
at  Aleksander Pe{\l}czy{\'n}ski  Memorial Conference
(13 July --  19 July, 2014,  B\c{e}dlewo, Poland)
as well as for very useful information which B. S. Mityagin gave me
in  letters. Also, I am very grateful to him
for drawing my attention to the paper of M. I. Zelikin [16]. %!!!REF

   %%%

\begin{thebibliography}{99}

  \bibitem{1} Gohberg I., Goldberg S., Krupnik N.,
{Traces and Determinants of Linear Operators},
Birkh\"auser Verlag, Basel-Boston-Berlin, 2000.

\bibitem{2}  Grothendieck A.,
 {Produits tensoriels topologiques et espaces nucl\'eaires},
  Mem. Amer. Math. Soc., 1955, 16.%, 196 + 140.

   \bibitem{3}  Grothendieck A.,
{La th\'eorie de Fredholm},
{Bull. Soc. Math. France},  1956, 84, 319--384.

 \bibitem{4}  K\"onig H.,
{Eigenvalues of Operators and Applications},
 Handbook of the geometry of Banach spaces, vol. 1, Chap. 22, 2001, 941--974.

 \bibitem{5}
 Lidski\v{\i} V. B.,
  {Nonselfadjoint operators having a trace},
 {Dokl. Akad. Nauk SSSR}, 1959, 125, 485--487.

 \bibitem{6}  Mityagin B. S.,
{Criterion for $Z_d$-symmetry of a Spectrum of a Compact Operator},
arXiv: math.FA/1504.05242.

 \bibitem{7} Mityagin B. S.,
{A criterion for the $\mathbb{Z}_d$-symmetry of the spectrum of a compact operator},
{J. Operator Theory}, 2016, 76, 57--65.

    \bibitem{8} Pietsch A.,
{Operator Ideals},
North Holland, 1980.

     \bibitem{9}  Pietsch A.,
{Eigenvalues and s-Numbers},
Cambridge studies in advanced mathematics,
{13}, Cambridge University Press, 1987.

 \bibitem{10} Reinov O. I.,
{A simple proof of two theorems of A. Grothendieck},
Vestn. Leningr. Univ. , 1983, 7, 115-116.

  \bibitem{11}
   Reinov O. I.,
     {Approximation properties $AP_s$ and p-nuclear operators
      (the case $0 <s\le 1)$},
 {J.  Math. Sciences}, 2003, 115, 2243--2250.

  \bibitem{12}  Reinov O. I.,
 {On linear operators with $s$-nuclear adjoints, $0<s\le1,$}
  {J. Math. Anal. Appl.}, 2014, 415,  816--824.

 \bibitem{13}
 Reinov O., Latif Q.,
 {Grothendieck-Lidski\v{\i} theorem for subspaces of quotients of  $L_p$-spaces},
 {Proceedings of the international mathematical conference Function Spaces X},
Banach Center Publications, 2014, 102, 189--195.

 \bibitem{14}  Schatten R.,
 {A theory of cross-spaces}, Princeton University Press, 1950.

 \bibitem{15} White M. C.,
{Analytic multivalued functions and spectral trace},
Math. Ann., 1996, 304, 665--683.

 \bibitem{16} Zelikin M. I.,
{A criterion for the symmetry of a spectrum},
{Dokl. Akad. Nauk}, 2008, 418, 737--740.

              \end{thebibliography}
  \end{document}